\documentclass[a4paper,12pt]{article}
\usepackage{amsthm,amsmath,amssymb,amsfonts}
\usepackage{float}
\usepackage{multirow}
\usepackage{multicol}

\theoremstyle{plain}

\theoremstyle{definition}
\newtheorem{definition}{Definition}
\newtheorem{example}{Example}

\begin{document}


\title{\bf Improving efficiency in fuzzy regression modeling by Stein-type shrinkage}

\author{{M. Kashani, M. Arashi\footnote{Corresponding author, Email: m\_arashi\_stat@yahoo.com} and M.R. Rabiei
}
\vspace{.5cm} \\\it
Department of Statistics, Faculty of Mathematical Sciences\\\it Shahrood University of Technology, Shahrood, Iran}

\date{}
\maketitle

\begin{abstract}
The fuzzy linear regression (FLR) modeling was first proposed making use of linear programming and then followed by many improvements in a variety of ways. In almost all approaches changing the meters, objective function, and restrictions caused to improve the fuzzy measure of efficiencies (FMEs). In this paper, from a totally different viewpoint, we apply shrinkage estimation strategy to improve FMEs in the FLR modeling. By several illustrative examples, we demonstrate the superiority of the proposed estimation method. In this respect, we show fuzzy shrinkage estimates improve FMEs estimation dramatically compared to the existing methods. 
\end{abstract}

\vspace{9pt} \noindent {\it Key words and phrases:} Fuzzy linear regression; Fuzzy least distance; Bootstrap; Shrinkage.
\par


	\section{Introduction}
The fuzzy linear regression (FLR) was introduced by \cite{1}, for the first time in 1982. This approach was an extension to the classical linear regression modeling and developed by several researchers in subsequent years from possibility view point \cite{a,b,c,e,f}. On the other hand, using least squares theory \cite{5} and then \cite{6} proposed another proposal of a different FLR modeling, called the fuzzy least squares regression. In 1991, Savic and Pedrich \cite{7} proposed the idea of two-stage methods, a combination method, in order to increase the efficiency of the Tanaka's approach by combining possibility and least squares view points. In the process of analyzing numerous problems in the study of FLR modeling and inspired by the methods used in classical regression, innovative methods such as robust methods \cite{8,9}, bootstrap resampling \cite{2,11} and etc. were also used to increase efficiency in FLR models. All of these approaches and perspectives were based on goodness of fit (GOF) criteria.

In this respect, the estimation of fuzzy regression parameters is of main importance, since the estimates will be used in the goodness of fit measures and prediction accuracy of the underlying model. Hence, improving estimation, in a direction of improving prediction or either goodness of fit, is very important. One approach of improving the estimation in this direction is to reduce the mean distance error (MDE) of estimation by shrinking.

Suppose the random vector $ \mathbf{X} $ of dimension $ p $, has $ p $-variate normal distribution with mean 
$ \boldsymbol{\mu} =(\mu_1,\mu_2,\dots, \mu_p)^\top\in \mathbb{R}^p $ 
and identity covariance matrix, denoted by
 $ \mathbf{X}\sim N_p (\boldsymbol{\mu}, \mathbf{I}_p) $.
 To estimate the location vector $ \boldsymbol{\mu} $, one may use 
 $ \hat{\boldsymbol{\mu}}= \boldsymbol{\delta} (\mathbf{X})=\mathbf{X} $.
 This estimator has many good properties. For it $ p<3 $, is the best linear unbiased estimator, minimax and admissible. However, \cite{4} showed when $ p\geq 3 $, it is inadmissible and in 1961, with his student, presented a more precise type of estimator known as James-Stein shrinkage estimator \cite{12}. The idea was to shrink $ \hat{\boldsymbol{\mu}} $ toward zero to obtain a minimax estimator.
 
Assume that the vector of estimators in a FLR model is 
 $ \hat{\tilde{\boldsymbol{\theta}}}=(\hat{\tilde{\theta}}_1, \hat{\tilde{\theta}}_2, \dots,\hat{\tilde{\theta}}_p )^\top $,
 then the fuzzy shrinkage estimator vector denoted by $ \hat{\tilde{\theta}} $ has form
 \begin{equation}
 \hat{\tilde{\boldsymbol{\theta}}}^{S} (k)=\left(\hat{\tilde{\theta}}_i^{S} (k),\, i=1,2,\dots,p \right)^\top,
 \end{equation}
 where 
 \begin{equation}
 \hat{\tilde{\theta}}_i^{S} (k)=\left(1-\dfrac{k}{\hat{\tilde{\theta}}_i^2}\right) \hat{\tilde{\theta}}_i,\qquad k>0,
 \end{equation} 
 $ k $ is called the shrinkage constant, which for a given value, the model error will be minimized. In this paper, we will be showing that with the use of this idea, with the least computational cost, the fuzzy GOF measures can be improved. Hence, we organize our paper as follows. In section \ref{sec:2} we define LR-fuzzy numbers and operations between them. Also explicit formula for the fuzzy Stein-type shrinkage estimator is given, while in section \ref{sec:3} the underlying fuzzy GOF measures are outlined. Extensive comparisons between the proposed method of estimation with some existing ones in the literature is conducted in section \ref{sec:4}. Our paper is concluded in section \ref{sec:5}. 
 \section{Definitions \& Concepts}\label{sec:2}
This section briefly reviews several concepts and terminologies related to fuzzy numbers and structure of JS. 
\subsection{Fuzzy preliminaries}
Let $ \chi $ be a universal set. A fuzzy set of $ \chi $ is a mapping
 $ \widetilde{A}:\chi \to [0,1] $, 
 which assigns a degree of membership $ 0 \leq \widetilde{A} (x)\leq 1 $ to each $ x\in \chi $. For each $ \alpha \in (0,1] $, the subset
 $ \{x\in \chi |\widetilde{A}(x)\geq \alpha \} $ 
 is called the $ \alpha $-cut of $ \widetilde{A} $ and is denote by $ \widetilde{A}_\alpha $. The set $ \widetilde{A}_0 $ is also defined as equal to the closure of
 $ \{x\in \mathbb{R} | \widetilde{A} (x )>0\} $.
 Let $ \mathbb{R} $ be the set of all real numbers. A fuzzy set $ \widetilde{A} $ of $ \mathbb{R} $ is called a fuzzy number if it satisfies the following two conditions: 
\begin{enumerate}
	\item 
	For each $ \alpha \in [0,1] $, the set $ \widetilde{A}_\alpha $ is a compact interval, which will be denoted by 
	$ \left[\widetilde{A}_\alpha^L, \widetilde{A}_\alpha^U\right] $. 
	Here, $ \widetilde{A}_\alpha^L=\inf \{x\in \mathbb{R} | \widetilde{A} (x )\geq \alpha \} $ and $ \widetilde{A}_\alpha^U=\sup \{x\in \mathbb{R} | \widetilde{A} (x )\geq \alpha \} $.
\item 
There is a unique real number $ x^*=x_{\widetilde{A}}^*\in \mathbb{R} $, such that $ \widetilde{A} (x^*)=1 $, i.e. $ \widetilde{A}^{-1}(1) $ is a singleton set.
\end{enumerate} 
The set of all fuzzy numbers with continuous membership functions is denoted by $ \mathcal{F}(\mathbb{R}) $. Notably, the most commonly used type of fuzzy numbers in $ \mathcal{F}(\mathbb{R}) $ is the so-called \textit{LR-fuzzy} numbers denoted by $ \widetilde{A}=(l_A,m_A,r_A )_{LR}, l_A,r_A>0 $, where $ m_A $, $ l_A $ and $ r_A $ are center, left spread and right spread of the LR-fuzzy number respectively.

The membership function of an LR-fuzzy number $ \widetilde{A} $ is defined by:
\begin{equation}
\widetilde{A}(x)=
\begin{cases}
L\left(\frac{m_A-x}{l_A}\right), & x\leq m_A,\\
R\left(\frac{x-m_A}{r_A} \right), & x>m_A,
\end{cases}
\end{equation}
	where $ L $ and $ R $ are strictly decreasing functions from $ [0,1] $ to $ [0,1] $ satisfying $ L(0)=R(0)=1 $ and $ L(1)=R(1)=0 $. For example special type of LR-fuzzy number is the triangular fuzzy number (TFN) with the shape functions 
	$ L(x )=R(x )=\max\{ 0, 1-|x|\}, x\in \mathbb{R} $, denoted by $ \widetilde{A}=(l_A, m_A, r_A)_T $. If $ l_A=r_A $, then $ \widetilde{A} $ is called a symmetric triangular fuzzy number. 
	
	Algebraic operations on fuzzy numbers that we use in this paper are defined based on the extension principle as follows.
\begin{definition}
		 \cite{13}. If $ \widetilde{A}=(l_A,m_A,r_A )_{LR} $ and 
		 $ \widetilde{B}=(l_B,m_B,r_B )_{LR} $ be two LR-type fuzzy numbers and 
		 $ \lambda \in \mathbb{R} $. Then: 
	\begin{enumerate}
		\item
		$ \widetilde{A}\oplus \widetilde{B}=(l_A+l_B,m_A+m_B,r_A+r_B )_{LR} $.
	\item
	$ \lambda \otimes \widetilde{A}=
	\begin{cases}
	(\lambda l_A,\lambda m_A,\lambda r_A )_{LR}, & \lambda >0, \\
	(-\lambda r_A,\lambda m_A,-\lambda l_A)_{LR}, & \lambda <0.
	\end{cases} $.
	\end{enumerate}
\end{definition}
For more details refer to \cite{13}. 
 \subsection{Fuzzy Stein-type shrinkage estimation}

		In statistical inference, one of the most important criteria for performance analysis of an estimator is the mean squared error (MSE). When estimators are unbiased, we will achieve the goal by minimizing the estimator variance, and in fact, the calculation will be simpler. However, there are non-market estimators that have the above-mentioned criterion less than the unbiased estimators \cite{14,15}. One of these estimators is the Stein-type shrinkage. For each of the fuzzy regression model, in this section, we define the Stein-type shrinkage estimator.
		
		Let $ \widetilde{Y}_i=(l_{Y_i},m_{Y_i},r_{Y_i} )_{LR}, $ for $ i=1,2,\dots, n $ be fuzzy observed responses and $ (X_{i1},X_{i2},\dots, X_{ip}), \, i=1,2,\dots,n $ are real independent variables from fuzzy regression model as follows:
		\begin{equation}
		\widetilde{Y}_i=\widetilde{A}_{0}+\widetilde{A}_{1} X_{i1}+\widetilde{A}_{2} X_{i2}+ \dots +\widetilde{A}_{p} X_{ip},
		\end{equation}
		where $ \widetilde{A}_j=(l_{A_j},m_{A_j},r_{A_j} )_{LR}, \, j=0, 1, \dots, p, $ are model coefficients. Now, if $ \hat{\tilde{A}}_j $ represent the fuzzy estimate of the coefficient $ \widetilde{A}_j $ in then, its Stein-type shrinkage estimate is given by
		\begin{equation}
		\hat{\tilde{A}}_j^{S}=\left(\hat{l}_{A_j}^{S}, \hat{m}_{A_j}^{S}, \hat{r}_{A_j}^{S} \right)_{LR},
		\end{equation}
		where
		\begin{equation}\label{equ:4}
		\begin{cases}
		 \hat{l}_{A_j}^{S}=\left(1-\dfrac{k}{\hat{l}_{A_j}^2} \right) \hat{l}_{A_j},\\[0.5cm]
		 \hat{m}_{A_j}^{S}=\left(1-\dfrac{k}{\widehat{m}_{A_j}^2} \right) \widehat{m}_{A_j},\\[0.5cm]
		 \hat{r}_{A_j}^{S}=\left(1-\dfrac{k}{\hat{r}_{A_j}^2} \right) \hat{r}_{A_j}.
		\end{cases}
		\end{equation}		
		and $ k>0 $ is the shrinkage coefficient (tuning parameter). An important point in the proposed estimator is that, with increasing the $ k $, initial estimator coefficient in Eq. \eqref{equ:4} may be negative. To fix this problem, Stein in (1966), defined the positive-rule Stein estimator. Here, we will define and use the fuzzy positive-rule Stein. Thus, for the spreads parameter, fuzzy positive-rule Stein-type shrinkage estimates have form 
\begin{align}\label{equ:5}
\begin{cases}
\hat{l}_{A_j}^{S+}=\left(1-\dfrac{k}{\hat{l}_{A_j}^{S+^2}} \right)^+ \hat{l}_{A_j}^{S+}=\max\left\{0,\left(1-\dfrac{k}{\hat{l}_{A_j}^{S+^2}} \right) \hat{l}_{A_j}^{S+} \right\},\\[0.5cm]
\hat{r}_{A_j}^{S+}=\left(1-\dfrac{k}{\hat{r}_{A_j}^{S+^2}} \right)^+ \hat{r}_{A_j}^{S+}=\max\left\{0,\left(1-\dfrac{k}{\hat{r}_{A_j}^{S+^2}} \right) \hat{r}_{A_j}^{S+} \right\}.
\end{cases}
\end{align}	
		 \section{Fuzzy GOF Measures}\label{sec:3}
			In this paper, some fuzzy measure of closeness are used for evaluating the performance of a FLR model. Let, $ \widetilde{Y}_i=(l_{Y_i},m_{Y_i},r_{Y_i} )_{LR} $ and $ \hat{\tilde{Y}}_i=(l_{\widehat{Y}_i },m_{\widehat{Y}_i },r_{\widehat{Y}_i } )_{LR} $ are $ i^{th} $ LR-fuzzy observation and estimation numbers, respectively. In sequel we give some well-known fuzzy GOF measures which we will be used in comparing the performance of our proposed method with others.
			
			I. Sadeghpour and Gien measure: \cite{16}\\
			\begin{equation}
			\begin{split}
			D_{p,q}&=\frac 1n \sum_{i=1}^nD_{p,q}^2 (i),\\
			&=\frac 1n \sum_{i=1}^n\left((1-q) \int_0^1\left(\left(\hat{\tilde{Y}}_i \right)_\alpha ^-- \left(\widetilde{Y}_i \right)_\alpha^- \right)^p d\alpha +q\int_0^1\left(\left(\hat{\tilde{Y}} \right)_\alpha^+- \left(\widetilde{Y}_i \right)_\alpha ^+ \right)^p d\alpha \right),
			\end{split}
			\end{equation}
			Where, $ \left(\widetilde{Y}_i \right)_\alpha ^+ $ and $ \left(\widetilde{Y}_i \right)_\alpha ^- $ the upper and lower limits of the observed response by the alpha-cut and also for the estimated response respectively. In the special case, for triangular fuzzy numbers $ (p=2, q=\frac 12)$ we have:
			 \begin{equation}
			\begin{split}
			D_{2,\frac 12}=\frac 1n \sum_{i=1}^n\Big(\frac 16 &\Big[\left(\left(m_{\widehat{Y}_i}-l_{\widehat{Y}_i}\right)-\left(m_{Y_i}-l_{Y_i}\right)\right)^2+2\left(m_{\widehat{Y}_i }-m_{Y_i} \right)^2\\
			&+\left(\left(m_{\widehat{Y}_i }+r_{\widehat{Y}_i }\right)-\left(m_{Y_i}+r_{Y_i}\right)\right)^2\\
			&+\left(\left(m_{\widehat{Y}_i }-l_{\widehat{Y}_i }\right)-\left(m_{Y_i}-l_{Y_i}\right)\right)\left(m_{\widehat{Y}_i }-m_{Y_i} \right)\\
			&+\left(\left(m_{\widehat{Y}_i }+r_{\widehat{Y}_i }\right)-\left(m_{Y_i}+r_{Y_i}\right)\right)\left(m_{\widehat{Y}_i }-m_{Y_i} \right)\Big]\Big),
			\end{split} 
			\end{equation}
			II. Hassanpour et al. measure: \cite{17}
			\begin{equation}
			D_H \left(\widetilde{Y}_i,\hat{\tilde{Y}}_i \right)=\left|m_{\widehat{Y}_i }-m_{Y_i} \right|+\left|r_{\widehat{Y}_i }-r_{Y_i} \right|+\left|l_{\widehat{Y}_i }-l_{Y_i} \right|,
			\end{equation}
			III. Kelkinnama and Taheri measure: \cite{18}
			\begin{equation}\label{equ:9}
			\begin{split}
			D_{LR} \left(\widetilde{Y}_i,\hat{\tilde{Y}}_i \right)=&\frac 13 \Big(\left|m_{\widehat{Y}_i }-m_{Y_i} \right|+\left|\left(m_{\widehat{Y}_i }+w_r r_{\widehat{Y}_i }\right)-\left(m_{Y_i}+w_r r_{Y_i}\right)\right|\\
			&+\left|\left(m_{\widehat{Y}_i }-w_l l_{\widehat{Y}_i }\right)-\left(m_{Y_i}-w_l l_{Y_i}\right)\right|\Big),
			\end{split}
			\end{equation}
			Where, $ w_r=\displaystyle\int_0^1R^{-1} (\alpha )d\alpha $, and $ w_l=\displaystyle\int_0^1L^{-1} (\alpha )d\alpha $. Which is for triangular fuzzy numbers:
		 \begin{equation}
			\begin{split}
			D_{LR} \left(\widetilde{Y}_i,\hat{\tilde{Y}}_i \right)=&\frac 13 \Big(\left|m_{\widehat{Y}_i }-m_{Y_i} \right|+\left|\left(m_{\widehat{Y}_i }+\frac 12 r_{\widehat{Y}_i }\right)-\left(m_{Y_i}+\frac 12 r_{Y_i}\right)\right|\\
			&+\left|\left(m_{\widehat{Y}_i }-\frac 12 l_{\widehat{Y}_i }\right)-\left(m_{Y_i}-\frac 12 l_{Y_i}\right)\right|\Big).
			\end{split}
			\end{equation}
		 \section{Numerical Demonstrations}\label{sec:4}
		 In this section, we provide some numerical illustrations to compare the performance of the proposed Stein-type shrinkage method with the existing ones in the literature. 
		 \begin{example}\label{exa:1}
		 	Table \ref{tab1} contains real data (Dataset 1), which was first used by \cite{22}. This Dataset, which is related to a study on the cognitive response time of a nuclear power plant control room crew to an abnormal event, consists of crisp input and symmetric triangular fuzzy output. This Dataset is also analysis by \cite{18}. Indeed, they illustrated their model smaller $ D_{LR} $ give by Eq. \eqref{equ:9} compared to \cite{25}. Here, we show if we use the Stein-type shrinkage estimate Eq. \eqref{equ:4} for centers and its positive-part in Eq. \eqref{equ:5} for spreads, then we obtain smaller $ D_{LR} $ values compared to \cite{18}, dominating all other abovementioned works.
		 	 
		 \end{example}
			\begin{table}[H]
				\caption{(Dataset 1) Here $ \hat{\tilde{Y}}_i $ and $ \hat{\tilde{Y}}_i^S $ represent the fitted values of 
					response $ \widetilde{Y}_i $ using \cite{18} and shrinkage methods, respectively.
				}\label{tab1}
				\begin{center}
					\begin{tabular}{ccccccc}
						\hline
			No.	& $ x_1 $ &	$ x_2 $ & $ x_3 $ &	$ \widetilde{Y}_i $ &$ \hat{\tilde{Y}}_i  $ &	$ \hat{\tilde{Y}}_i^S  $\\
			\hline
			1&	2.00&	0.00&	15.25&	$ (5.83,3.56)_T $&	$ (6.97,1.78)_T $&	$ (6.96,0.94)_T $\\
			2&	0.00&	5.00&	14.13&	$ (0.85,0.52)_T $&	$ (0.85,2.29)_T $&	$ (0.83,1.45)_T $\\
			3&	1.13&	1.50&	14.13&	$ (13.93,8.50)_T $&	$ (4.05,1.80)_T $&	$ (4.04,1.00)_T $\\
			4&	2.00&	1.25&	13.63&	$ (4.00,2.44)_T $&	$ (3.35,1.92)_T $&	$ (3.35,1.14)_T $\\
			5&	2.19&	3.75&	14.75&	$ (1.58,0.96)_T $&	$ (2.46,2.61)_T $&	$ (2.47,1.72)_T $\\
			6&	0.25&	3.50&	13.75&	$ (1.58,0.96)_T $&	$ (1.72,1.99)_T $&	$ (1.70,1.19)_T $\\
			7&	0.75&	5.25&	15.25&	$ (8.18,4.99)_T $&	$ (2.06,2.63)_T $&	$ (2.05,1.71)_T $\\
			8&	4.25&	2.00&	13.50&	$ (1.85,1.13)_T $&	$ (1.85,2.64)_T $&	$ (1.89,1.81)_T $\\
			\hline
					\end{tabular}
				\end{center}
			\end{table}
			Using \cite{18} and equations \eqref{equ:4} and \eqref{equ:5}, the fuzzy regression models are given by \eqref{(11a)} and \eqref{(11b)}. 
		 Here, $ \oplus $ and $ \odot $ are the specific fuzzy sum and product definitions used in \cite{18}. The fuzzy GOF measures are tabulated in Table \ref{tab2}.
			\begin{align}
			\hat{\tilde{Y}}_i=(&-14.8998,0.2500)_T\oplus(-0.2505, 0.2500)_T \odot x_1 \notag \\
			&\oplus(-0.9558,0.2216)_T \odot x_2\oplus(1.4670,0.0837)_T \odot x_3 \tag{13a}\label{(11a)}\\
			\hat{\tilde{Y}}_i^S=(&-14.8995,0.2324)_T\oplus(-0.2329, 0.2324)_T \odot x_1 \notag\\
			&\oplus(-0.99137,0.2017)_T \odot x_2\oplus(1.4640,0.0310)_T \odot x_3 \tag{13b}\label{(11b)}
			\end{align}
				\begin{table}[H]
					\caption{Comparison of the shrinkage method for $ k=0.0044 $ with \cite{18} using the $ D_{LR} $ as the fuzzy GOF measure.}\label{tab2}
					\begin{center}
						\begin{tabular}{l@{\hspace*{8mm}}l}
							\hline
							Model&	$ D_{LR} $\\
							\hline
							Kelkinnama and Taheri \cite{18}&20.1521\\
							Shrinkage method&	19.4929\\							
							\hline
						\end{tabular}
					\end{center}
				\end{table}	
			According to the Table \ref{tab2}, result of the $ D_{LR} $ for the shrinkage method is smaller than \cite{18}. In addition for the shrinkage constant in the boundary $ (0,0.0308] $, termed as the optimal boundary, the proposed method is still superior. Note that $ k=0.0044 $ is the value for which $ D_{LR} $ has the smallest value.
				
				\begin{example}\label{ex:2}
					Here, in Table \ref{tab3} (Dataset 2), the data used by Tanaka et al. \cite{24}.
				\end{example} 
				In this example, we want to demonstrate the effects of the suggested methodology on innovative techniques. Therefore, using the Bootstrap method in least square view, we obtain the estimated responses and corresponding contraction values for them under the suitable k-value. Then we compare the results of two approach in Table \ref{tab4}.
				
				\begin{table}[H]
					\caption{(Dataset 2) Here $ \hat{\tilde{Y}}_i $ and $ \hat{\tilde{Y}}_i^S $ represent the fitted values of 
						response $ \widetilde{Y}_i $ using Bootstrap and shrinkage methods, respectively.
					}\label{tab3}
					\begin{center}
						\begin{tabular}{ccccc}
							\hline
							No. &	$ X_i $	& $ \widetilde{Y}_i $&	$ \hat{\tilde{Y}}_i $&	$ \hat{\tilde{Y}}_i^{S}  $\\
							\hline
							1&	1&	$ (8.00,1.80)_T $&	$ (6.72,2.00)_T $&
							$ (6.65,1.79)_T $
\\
							2&	2&	$ (6.40,2.20)_T $&	$ (8.41,2.16)_T $&
							$ (8.27,1.79)_T $
\\
							3&	3&	$ (9.50,2.60)_T $&	$ (10.09,2.32)_T $&	$ (9.90,1.79)_T $
\\
							4&	4&	$ (13.50,2.60)_T $&	$ (11.78,2.47)_T $&	$ (11.53,1.79)_T $
\\
							5&	5&	$ (13.00,2.40)_T $&	$ (13.47,2.63)_T $&	$ (13.16,1.79)_T $
\\
							\hline
						\end{tabular}
					\end{center}
				\end{table}
In this example, using the dataset 2, we want to show the effect of improving the our method versus the Boot technique. So, the models derived from these two techniques will be as follows
\begin{align}
\hat{\tilde{Y}}_i=(5.0365,1.8469)_T+(1.6862, 0.1565)_T X_i \tag{14a}\\
\hat{\tilde{Y}}_i^{S}=(5.0172,1.7943)_T+(1.6285, 0.000)_T X_i \tag{14b}
\end{align}			
\begin{table}[H]
	\caption{Comparison of the shrinkage method for $ k=0.0972 $ with \cite{23} using the $ D_{LR} $ as the fuzzy GOF measure.}\label{tab4}
	\begin{center}
		\begin{tabular}{l@{\hspace*{9mm}}l}
			\hline
			Model &	$ D_{LR} $ \\
			\hline
			Arabpour and Moradi \cite{2}&	6.06747\\
			Shrinkage method&	5.85522\\			
			\hline
		\end{tabular}
	\end{center}
\end{table}	
According to the Table \ref{tab3}, result of the $ D_{LR} $ for the shrinkage method is smaller than \cite{2} and addition for the shrinkage constant in the boundary $ (0,0.2138] $. The optimal boundary, the proposed method is still superior. Note that $ k=0.0972 $ is the value for which $ D_{LR} $ has the smallest value.

As can be seen from the results of Table \ref{tab4}, while the boot method is one of the most powerful methods for estimating regression model coefficients in low-volume, high-performance samples, but the use of JS method leads to a decrease in $ D_{LR} $ criteria and efficiency The model rises.
\begin{example}
 Consider the cheese data in Table \ref{tab5} (Dataset 3). This data is presented based on the quality of cheese tasting as a response variable, which is evaluated by a specialist. To generate the fuzzy triangular response, we add 15\% of the center points to the spreads. The explanatory variables are acetic acid (Acetic), sulfuric acid ($ H_2S $), and lactic acid (Lactic).
\end{example}
\begin{table}[H]
	\caption{(Dataset 3) Here $ \hat{\tilde{Y}}_i $ and $ \hat{\tilde{Y}}_i^S $ represent the fitted values of 
		response $ \widetilde{Y}_i $ using \cite{21} and shrinkage methods, respectively.}\label{tab5}
	\begin{center}
			\begin{tabular}{ccccccc}
				\hline
				No.&	Acetic($ x_1 $)&	$ H_2S(x_2) $	&Lactic($ x_3 $)&	$ \widetilde{Y}_i $ &	$ \hat{\tilde{Y}}_i $&	$ \hat{\tilde{Y}}_i^S $\\
				\hline
				1&	4.543&	3.135&	0.86&	$ (12.30,1.845)_T $ &	$ (6.89,1.243)_T $ &
				$ (7.49,1.390)_T $ \\
				2&	5.159&	5.043&	1.53&	$ (20.90,3.135)_T $ &	$ (22.34,2.059)_T $ &
				$ (23.32,2.304)_T $ \\
				3&	5.366&	5.438&	1.57&	$ (39.00,5.850)_T $ &	$ (27.74,2.188)_T $ &
				$ (28.83,2.447)_T $ \\
				4&	5.759&	7.496&	1.81&	$ (47.90,7.185)_T $ &	$ (34.62,2.874)_T $ &
				$ (36.30,3.206)_T $ \\
				5&	4.663&	3.807&	0.99&	$ (5.60,0.840)_T $ &	$ (9.02,1.488)_T $ &
				$ (9.80,1.662)_T $ \\
				6&	5.697&	7.601&	1.09&	$ (25.90,3.885)_T $ &	$ (30.40,2.616)_T $ &
				$ (32.36,2.900)_T $ \\
				7&	5.892&	8.726&	1.29&	$ (37.30,5.595)_T $ &	$ (33.73,3.018)_T $ &
				$ (35.99,3.347)_T $ \\
				8&	6.078&	7.966&	1.78&	$ (21.90,3.285)_T $ &	$ (43.09,2.997)_T $ &
				$ (44.92,3.340)_T $ \\
				9&	4.898&	3.85&	1.29&	$ (18.10,2.715)_T $ &	$ (17.04,1.621)_T $ &
				$ (17.72,1.817)_T $ \\
				10&	5.242&	4.176&	1.58&	$ (21.00,3.150)_T $ &	$ (27.59,1.830)_T $ &
				$ (28.26,2.056)_T $ \\
				11&	5.74&	6.142&	1.68&	$ (34.90,5.235)_T $ &	$ (37.62,2.434)_T $ &
				$ (38.90,2.720)_T $ \\
				12&	6.446&	7.908&	1.9&	$ (57.20,8.580)_T $ &	$ (55.04,3.028)_T $ &
				$ (56.80,3.378)_T $ \\
				13&	4.477&	2.996&	1.06&	$ (0.70,0.105)_T $ &	$ (5.79,1.279)_T $ &	$ (6.284,1.440)_T $ \\
				14&	5.236&	4.942&	1.3	& $ (25.90,3.885)_T $ &	$ (24.39,1.937)_T $ &
				$ (25.42,2.164)_T $ \\
				15&	6.151&	6.752&	1.52&	$ (54.90,8.235)_T $ &	$ (48.19,2.545)_T $ &
				$ (49.71,2.836)_T $ 
				\\
				\hline
		\end{tabular}
	\end{center}
\end{table}
Here, we follow the method of \cite{21}. The fitted fuzzy regression models based on least squares and shrinkage methods are given by \eqref{(13a)} and \eqref{(13b)}, respectively.
				\begin{align}
				\hat{\tilde{Y}}_i=(&0,-127.6929,0)_T+(0,31.1153, 0)_T x_1\notag\\
				&+(0.57328,-2.9192,0)_T x_2+(0.8013,2.7644,0)_T x_3\tag{15a}\label{(13a)}\\
				\hat{\tilde{Y}}^S_i=(&0,-127.6854,0)_T+(0,31.0843, 0)_T x_1\notag\\
				&+(0.6276,-2.5886,0)_T x_2+(0.9438,2.7644,0)_T x_3 \tag{15b}\label{(13b)}
				\end{align}
	Result of comparison between the shrinkage method with \cite{21} are summarized in Table \ref{tab6}. Here, we used the fuzzy GOF measures of \cite{16,17,18}. 
\begin{table}[H]
	\caption{Comparison of the shrinkage method with \cite{21} using the fuzzy GOF measure.}\label{tab6}
	\begin{center}
			\begin{tabular}{lccc}
				\hline 
				Model&	\multicolumn{3}{c}{GOF measures and values}\\
				\hline
				\multirow{2}{*}{Hassanpour et al. \cite{21}}& $ D_{LR} $&	$ D_{2, \frac 12} $	& $ D_H $\\
				&	89.9129&	68.3101	&157.9474\\
				\hline
				\multirow{2}{*}{Shrinkage}&	88.0382&65.0767&146.2433\\
			 &$ (k=1.183) $	&$ (k=0.965) $	& $ (k=1.524) $\\
				Optimal boundary&$ (0,1.759] $	&$ (0,1.929] $&	$ (0,4.335] $
\\				
				\hline
		\end{tabular}
	\end{center}
\end{table}	
According to Table \ref{tab6}, result of value of the GOF measures for all show that the shrinkage method is superior to \cite{21}. 
\begin{example}\label{exa:4}
	 In this example (Dataset 4), we will examine another form of FLR models. The data is taken from \cite{a} and the FLR moder has form
\begin{equation}\label{equ:10}
\widetilde{Y}_i=\widetilde{A}_0+\widetilde{A}_1 \widetilde{X}_{i1}+\widetilde{A}_2 \widetilde{X}_{i2}+ \dots+\widetilde{A}_P \widetilde{X}_{iP}, \tag{16}
\end{equation} 
where $ \widetilde{Y}_i=(l_{Y_i}, m_{Y_i}, r_{Y_i})_T $ and $ \widetilde{X}_i=(l_{X_i}, m_{X_i}, r_{X_i})_T $ for $ i=1, 2, \dots,n $ are fuzzy output and input observations.
\end{example}
\begin{table}[H]
	\caption{(Dataset 4) Here $ \hat{\tilde{Y}}_i $ and $ \hat{\tilde{Y}}_i^S $ represent the fitted values of response $ \widetilde{Y}_i $ using \cite{f} and shrinkage method, respectively.}\label{tab7}
\begin{center}
	\begin{tabular}{ccccc}
		\hline
No.&	$ \widetilde{X}_i $& $ \widetilde{Y}_i $	&	$ \hat{\tilde{Y}}_i $&	$ \hat{\tilde{Y}}_i^S $\\
\hline
1&	$ (2.00,0.50)_T $ & $ (2.00,0.50)_T $ & $ (4.68,0.50)_T $& 	$ (4.52,0.48)_T $ \\
2&	$ (3.50,0.50)_T $ & $ (3.50,0.50)_T $ & $ (5.51,0.50)_T $ &	$ (5.23,0.48)_T $ 
\\
3&	$ (5.50,1.00)_T $ &	$ (5.50,1.00)_T $ &	$ (6.61,1.00)_T $ &	$ (6.18,0.96)_T $ \\
4&	$ (7.00,0.50)_T $ & $ (7.00,0.50)_T $ & $ (7.43,0.50)_T $ &	$ (6.90,0.48)_T $ \\
5&	$ (8.50,0.50)_T $ &	$ (8.50,0.50)_T $ &	$ (8.26,0.50)_T $ &	$ (7.61,0.48)_T $ \\
6&	$ (10.50,1.00)_T $ & $ (10.50,1.00)_T $& 	$ (9.36,1.00)_T $& 	$ (8.56,0.96)_T $ \\
7&	$ (11.00,0.50)_T $ & $ (11.00,0.50)_T $ &	$ (9.63,0.50)_T $ &	$ (8.80,0.48)_T $ \\
8&	$ (12.50,0.50)_T $ &$ (12.50,0.50)_T $ 	& $ (10.46,0.50)_T $ &	$ (9.51,0.48)_T $ \\
\hline
	\end{tabular}
\end{center}
\end{table}

This data has been examined and modeled by \cite{f}. In order to distinguish between methods, we present the model estimated by \cite{f} with $ \hat{\tilde{Y}}_i $, and the corresponding contraction model with $ \hat{\tilde{Y}}_i^S $. 
\begin{align}
\hat{\tilde{Y}}_i= (3.58,0.00)_T + (0.55, 1.00)_T \widetilde{X}_i\tag{17a}\\
\hat{\tilde{Y}}_i^S= (3.57,0.00)_T + (0.48, 0.96)_T \widetilde{X}_i\tag{17b}	
\end{align}
Result of comparison between the shrinkage method with \cite{f} are summarized in Table \ref{tab8}. Here, we used the fuzzy GOF measures of \cite{16,17,18}. 
	\begin{table}[H]
		\caption{Comparison of the shrinkage method with \cite{f} using the fuzzy GOF measure.}\label{tab8}
		\begin{center}
		\begin{tabular}{lccc}
			\hline 
			Model&	\multicolumn{3}{c}{GOF measures and values}\\
			\hline
			\multirow{2}{*}{Nasrabadi and Nasrabadi \cite{f}}& $ D_{LR} $&	$ D_{2, \frac 12} $	& $ D_H $\\
			&	6.9350&	5.6933&	7.6550\\
			\hline
			\multirow{2}{*}{Shrinkage}&	5.6640&	5.1435&6.2759\\
			&$ (k=0.041) $&	$ (k=0.017) $&$ (k=0.062) $\\
			Optimal boundary&$ (0,0.048] $&	$ (0,0.034] $ &	$ (0,0.092] $

\\				
			\hline
		\end{tabular}
		\end{center}
	\end{table}
According to the results of Table \ref{tab8}, our shrinkage approach is superior. 

The results of Table \ref{tab8} show that using the our method, all of our contiguity criteria for two fuzzy numbers, with high accuracy in \cite{f}, have been reduced. It should be noted that the positive effects of this method are calculated with a small cost. In addition, this method allows a specialist to have a set of points. Reduced to two fuzzy numbers with high accuracy in shrinkage method. It should be noted that the positive effects of this method are calculated with a low-cost. In addition, this method allows a specialist to have a set of points.

\section{Conclusions}\label{sec:5}
In the present work, we incorporated a prior knowledge embedded in the shrinkage strategy in the FLR modeling to improve the GOF measures. In other words, we allowed the fuzzy estimates (obtained from any fuzzy method) be shrunken toward the origin and specifically we used the Stein-type shrinkage estimator for the center points and its positive part for spreads of the fuzzy estimates. The reason we used the positive part estimator is the spreads must be positive.

According to illustrative examples, the proposed method works will in the sense of providing smaller distance measures. The computational cost of this method is fairly small and is of the same order as original fuzzy estimates.

Another benefit of the proposed approach is the optimal boundary for the shrinkage constant. As we shrink the spreads they become smaller, however, specialist may desire to have spreads with some actual size. Optimal boundaries allow the specialist to select a shrinkage constant which still serves the superiority, but with satisfactory spread. 

For future directions, first we point that the proposed method is very flexible and allows the specialist to use the shrinkage method on each of the center, left spread or right spread according to the usage, separately or simultaneously. Hence, for further work, one may extend our results for the case where the specialist wants to shrink the estimates towards an initial guess, for practical purposes. Further, our proposed method can be extended for LR fuzzy number and FLR models with error term (see \cite{27,28} for details).



\end{document}